
\input amstex

\documentstyle{amsppt}
\magnification =1200
\baselineskip=18pt
\NoBlackBoxes
\def\ls{\vskip.25in}
\def\ss{\vskip.15in}
\def\O{\Cal O}
\def\C{\Bbb C}
\def\H{\Bbb H}
\def\g{\goth g}
\def\R{\Cal R}
\def\P{\Bbb P}
\def\lra{\longrightarrow}

\def\Spec{\text{\rm Spec}}

\def\gr{\text{\rm gr}}
\def\im{\text{\rm im}}
\def\ker{\text{\rm ker}}
\def\X#1{X\kern-.5em<\kern-.3em#1\kern-.3em>}
\def\opp{\operatornamewithlimits{\oplus}}

\centerline{\bf How to tell you're hearing a Calabi-Yau:}
\centerline{\bf Universal variations of Hodge structure and local}
\centerline{\bf Schottky relations for Calabi-Yau manifolds}
\vskip.15in
\centerline{\bf Z. Ran}
\ls
The purpose of this paper is twofold.  First, we give a canonical formula for
the variation of Hodge structure associated to the $m$-th order universal
deformation of
an arbitrary compact K\"ahler manifold, this variation being viewed as a module
over the
base of the deformation.  Second, we specialize to the case of a Calabi-Yau
manifold
$X$ where we give a formula for the $m$-th differential of the period map of
$X$ and deduce formal defining equations for its image (Schottky relations);
these
are (necessarily infinite, in dimension $\geq 3$) power series in the middle
cohomology.

We will use the method of canonical infinitesimal deformations, developed by
the author
in earlier papers [R1, R2].  This method gives a canonical description of
infinitesimal
moduli spaces and, what's more, natural maps involving them.  While it might be
argued that a germ of a smooth space--such as the moduli of an unobstructed
manifold--is
a rather rigid featureless object, making a canonical description of it
uninteresting, on the contrary {\it maps} involving such germs can be quite
interesting; in the case of
moduli, the method of canonical infinitesimal deformations provides a vehicle
for studying
such maps.  For instance in the case at hand the $n$-th derivative of the
period
map of a Calabi-Yau $n$-fold $X$ is a filtered map
$$
T^n_X M \to H^n_{DR} (X) /F^n
$$
whose associated gradeds $S^i H^1(T_X) \to H^{n-i,i}_X$ are the so-called
Yukawa-Green
forms (cf. [G]).  We will develop cohomological formulas for this and other
derivative
maps (Theorem 3 below), which will allow us to determine their image and derive
(Schottky) relations defining this image, essentially in terms of some
generalized
Yukawa-Green type forms (corollary 3.1 below).  For $n = 2$ we recover the
celebrated
`period quadric' of $K3$ surface theory; for $n\geq 3$ the relations seem to be
new.
For $n=3$ the situation is particularly interesting assuming the 'mirror
conjecture'
because then the higher derivatives of the period map, here computed, are
related- in fact,
carry equivalent information to- the 'quantum cohomology' (esp. numbers of
rational
curves, etc.) of the mirror of $X$. We hope to return to this in greater detail
elsewhere.

The present methods should be applicable in other Schottky-type problems:  the
case of
curves is being developed by
G. Liu (UCR dissertation, to appear).

This paper is a revised version of a manuscript entitled `Linear structure on
Calabi-Yau
moduli spaces' (May 1993).  We are grateful to Professors P. Deligne and M.
Green for
their enlightening comments .

\subheading{1.  Preliminaries}
\ss
\noindent
{\it 1.1  Functors on $S$--Modules}
\medskip
In [R1] we showed how an Artin local $\Bbb C$--algebra may be reconstructed
from a certain `order--symbolic' or OS structure on the space of
$\Bbb C$--valued differential operators on $S$.
Our purpose here is to note an analogue of this for $S$--modules.

Now fix a local $\Bbb C$--algebra $S$ with maximal ideal $\goth m$ and
residue field $S/\goth m = \Bbb C$ and put
$$
B_0^i = D^i(S,\Bbb C) = \text{\rm Hom}(S_i,\Bbb C) = S_i^*,
\kern2em S_i = S/\goth m^{i+1},
$$
and
$$
B^i = D_+^i(S,\Bbb C) = (\goth m/\goth m^{i+1})^*.
$$
For an $S$--module $E$, put
$$
B^i(E) = B_0^i\otimes_S E,
$$
where $B^i_0$ is viewed as $S$--bi--module and $E$ as (symmetric)
$S$--bi--module.  At least when $E\otimes S_i$ is $S_i$--free,
$B^i(E)$ may be identified with the right $S$--module of differential
operators $D^i(E^\vee,\Bbb C)$, $E^\vee=\text{\rm Hom}_S(E,S)$.

We have a symbol map
$$
\sigma^i: B_0^i \rightarrow B^i\otimes_{\Bbb C}B^{i-1}_0
$$
which factors through $F_i(B^i\otimes_{\Bbb C}B^{i-1})$, where
$F_i$ is the filtration induced by the order filtration on $B^i$,
and this gives rise to a symbol map
$$
\sigma^i_E: B^i(E) \rightarrow B^i\otimes_{\Bbb C}B^{i-1}(E),
$$
which again factors through $F_i(B^i\otimes_{\Bbb C}B^{i-1}(E))$.
These $\sigma^i_E$, $i\le m$, together with the obvious maps
$B^0(E) \rightarrow B^1(E)\rightarrow \ldots\rightarrow B^m(E)$
are referred to as a `modular order--symbolic' (MOS) structure on
$B^m(E)$.  Note that $B^m(E)$ itself is a right $S_m$--module,
called the $m$--th transpose of $E$.

``Dually,'' suppose we are given an MOS structure $G^\cdot$, $G^i$ a right
$S_i$--module.  We then define an $S_m$--module $C^m(G^\cdot)$, called
the module of quasi--scalar homomorphisms $B^m_0 \rightarrow C^m(G^\cdot)$,
inductively as follows.
$$
C^0(G) = G^0
$$
$$
C^i(G) = \text{\rm all right~} S\text{\rm --linear maps~}
\varphi^i: B^i_0 \rightarrow G^i
$$
such that for some $\varphi^{i-1}: B^{i-1}_0\rightarrow G^{i-1}$ the following
diagrams commute.
$$
\CD
B_0^{i-1} @>{\varphi^{i-1}}>> G^{i-1} &&&&
\kern10em&& B_0^i @>{\varphi^i}>> G^i \\
@VVV @VVV &&
\kern10em&& @V{\sigma^i}VV @VV{\sigma^i_G}V  \\
B_0^i @>{\varphi^i}>> G^i &&&&
\kern10em&& B^i\otimes B^{i-1}_0 @>{\text{\rm id}\otimes\varphi^{i-1}}>>
B^i\otimes G^{i-1} \\
\endCD
$$

Note that we have natural maps
$$
E\rightarrow C^m(B^m(E))
\tag1.1
$$
$$
B^m(C^m(G^\cdot)) \rightarrow G^\cdot.
\tag1.2
$$
At least when $E$ is $S_m$--free (resp. $G^\cdot$ is `co--free',
i.e. a sum of copies of $B^m_0$ with the standard MOS structure),
these are isomorphisms.

\subsubhead
{1.2 Derivatives}
\endsubsubhead
\ss
For later use we want to give a more `geometric' interpretation of
$B^i (E)$, assuming $E$ corresponds to a geometric vector bundle $V(E)$
over a pointed space $(M,0)$.  The $(i+1)$st tangent sapce $T_{0,0}^{i+1}
(V(E))$ decomposes into components according to vertical degree ($=$
homogeneity
degree with respect to the natural $\C^*$ action), and $B^i(E)$ is just the
component
of vertical degree 1.  This may be verified easily.  Now suppose given a
`section', i.e.
an element $e \in E$, corresponding to a geometric cross-section $s_e : M \to
V(E)$ with
value $e(0) = s_e (0) \in E(0)$, and also to a map $u: M\otimes \C \to V(E) ,
u(x,t) =
(x, ts_e )$.  Then the natural map $B^i (e): B^i \to B^i (E)$ may be identified
as the degree 1 component of the composite
$$
T_0^i M {\overset x\partial/\partial t \to \longrightarrow} T_{(0,0)}^{i+1} (M
\times
\C) {\overset d_{(0,0)}^{i+1} u\to \longrightarrow} T_{(0,0)}^{i+1}
(V(E));
$$
or, what is the same under the natural identification
$
T_{(0,)}^{i+1} (V(E)) {\overset \sim\to\rightarrow} T_{(0,e(0))}^{i+1} (V(E))$
given by
translation by $S_e, B^i(e)$ is the vertical degree-1 component of the
composite
$$
T_0^i M {\overset d^is_e\to\longrightarrow} T_{0,e(0)}^i (V(E)) {\overset
Rx\to\rightarrow}
T_{0,e(0)}^{i+1} (V(E))
$$
where $R$ is the `Euler' vector field, corresponding to the $\C^*$-action.
Given a
trivialisation $V(E) {\overset \sim \to \rightarrow} E(0) \times M$,
corresponding to a
splitting $\sigma : B^i (E) \to E(0)$, so that $s_e$ is the graph of a function
$f: M \to E$, the composite $\sigma. B^i (e): B^i \to E(0)$ clearly coincides
with the
$i$-th differential of the function $f$, i.e. the map $d^i f: B^i = T^i M
{\overset T^i (f)\to \longrightarrow} T^i (E(0)) = \opp\limits_{1}^{i} S^i
(E(0)) \to
E(0)$.

\medskip
\subsubhead
{1.3.  Jacobi Complexes and Universal Deformations.}
\endsubsubhead
\ss
Fix a base space $X$ which for convenience we assume to be a compact
complex space (although the construction works more generally), and a
sheaf of Lie algebras on $X$, such that $H^0(\g)=0$.
As in [R1], we have Jacobi complexes $J^\cdot_m(\g)$ which may be
described as follows.  Let $\X m$ be the $m$--fold very
symmetric product of $X$, i.e. the space of subsets of $X$ of
cardinality $\in [1,m]$, with the topology induced by the natural map
$X^m \rightarrow \X m$.  For $i\le m$ let $\lambda^i(\g)$ be
the image of the exterior alternating product of $\g$, supported on
$\X i \subset \X m$.  Then the bracket on $\g$ gives
rise to a map $\lambda^i(\g)\rightarrow \lambda^{i-1}(\g)$,
and these fit together to form the complex $J^\cdot_m(\g)$---
$$
\lambda^m(\g)\rightarrow\lambda^{m-1}(\g)\rightarrow
\ldots\rightarrow\lambda^1(\g) = \g,
$$
in which we put $\lambda^i(\g)$ in degree $-i$.  Similarly, the
action of $\g$ on $E$ gives rise to a complex
$J^\cdot_m(\g,E)$ on $\X m\times X$---
$$
\lambda^m(\g)\boxtimes E\rightarrow\ldots\rightarrow
\g\boxtimes E \rightarrow E,
$$
in degrees $\in [-m,0]$, where the last term $E$ is supported on the
diagonal in $\X1\times X = X\times X$.  The natural maps
$$
J_i^\cdot(\g)\rightarrow J^\cdot_m(\g),\kern5em i\le m,
$$
$$
J^\cdot_m(\g)\rightarrow F_m(J^\cdot_{m-1}(\g)\boxtimes
J^\cdot_{m-1}(\g))
$$
give rise to an OS structure on $V_m = \Bbb H^0(J^\cdot_m(\g))$ which,
as in [R1], yields a $\Bbb C$--algebra structure on
$R_m=\Bbb C\oplus V_m^*$.
Similarly,
$$
G^m=\Bbb R^0p_{2*}(J^\cdot_m(\g,E)).
$$
Forms a sheaf of $R_m$--modules, and the natural map
$$
J^\cdot_m(\g,E)\rightarrow F_m(J^\cdot_m(\g)\boxtimes
J^\cdot_{m-1}(\g,E))
$$
endows $G^m$ with a MOS structure compatible with the OS structure
on $V^m$, whence as in $\S$ 1.2 a $R_m$-module
$$
C^m = C^m(G^m).
$$

In case $E$ itself is a $\C$-algebra on which $\g$ acts by derivations (or more
generally a
graded $\C$-algebra on which $\g$  acts by graded derivations), we have a
multiplication
map
$$
S^2 E \to E
$$
letting $\g$ act in the obvious way on $S^2E$, this map is clearly $\g$-linear,
hence
extends to a map of Jacobi complexes
$$
J^._m (\g, S^2 E) \to J_m^\cdot (\g, E).
$$
Combined with the natural map $\sigma^2 (J_m^\cdot (\g, E) \to J_m^\cdot (\g,
S^2 E),$
this gives
rise to a multiplicative structure on $R^0 p_2^* J_m^\cdot (\g, E)$, whence a
structure of
$R_m$-algebra on $C^m$.

Note also that $\g \oplus E$ has a natural structure of (differential graded)
Lie algebra,
and the image of $J_m^\cdot (\g, E)$ on $X < m+1>$ may be identified with a
direct
summand of $J_{m+1}^\cdot (\g \oplus E)$.

\ls
\subheading{2.  Universal variations of Hodge structure}
\ss
Fix a compact complex manifold $X$ with holomorphic tangent sheaf $\Theta =
\Theta_X$.
Note that $\Theta$ acts or $\O_X$ and more generally on the holomorphic De Rham
complex
$\Omega_X^\cdot$ by Lie derivative, so that the above constructions are
applicable.  Our
purpose is to apply them to give a description of the universal $m$-th order
variation of
Hodge structure of $X$, viewed as a module over the base of the $m$-th order
universal deformation of $X$.  We begin with a description of this deformation,
essentially
implicit in [R1].

\proclaim{Theorem 1}  Assuming $H^0 (\Theta) = 0$, the universal $m$-th order
deformation $X_m/R_m$ of $X$ and its relative DeRham complex
$\Omega^\cdot_{X_m/R_m}$ are
given by
\endproclaim
$$
\align
R_m &= \ \ \text{algebra associated to the}\quad OS \tag2.1\\
&\quad \text{structure on}\ \ \H^0 (J_m^{\cdot}(\Theta))\\
\Omega^\cdot_{X_m/R_m} &=\quad C_{R_m}^{m} (\ \ \R^{0,\cdot} p_{2*}
(J_m^{\cdot\cdot} (\Theta,
\Omega_X^\cdot ))\tag2.2
\endalign
$$
where $\Bbb R^{0,\cdot}$ means hyperderived image in the vertical $(J^\cdot)$
variable.

\demo{Proof}  (2.1) is proven in [R1], where a construction for $\O_{X_m}$ is
given, and it is
straightforward matter to check that this coincides with (2.2) for $\cdot =0$.
For
$\cdot = 1$, consider the $\Theta$-linear derivation $d: \O_X \to \Omega^1_X$.
By functoriality this yields an $R_m$-linear derivation $\O_{X_m} \to C^m_{R_m}
(R^0_{p_{2*}} (J_m^\cdot
\Theta, \Omega_X^1))$, hence an $\O_{X_m}$-linear map $\Omega_{X_m/R_m}^\cdot
\to
C^m_{R_m} (R^0_{p_{2*}} (J_m^\cdot (\Theta, \Omega_X^1))$, and by checking on
gradeds with
respect to the $m$-adic filtration on $R_m$ we see that this is an isomorphism.
 The
case of the rest of $\Omega^\cdot$ now follows easily.

Now we take up variations of Hodge structure.  The $R_m$-module $H^r_{DR}
(X_m/R_m) :=
\H^r(X, \Omega^\cdot_{X_m/R_m})$ is endowed with a (Hodge) filtration
$F^\cdot$, induced by
the stupid filtration on $\Omega^\cdot_{X_m/R_m}$.  At least when $X$ is
K\"ahler, so that
$H^r_{DR} (X_m/R_m)$ is $R_m$-free by Deligne [D], the pair $(H^r_{DR}
(X_m/R_m), F^\cdot)$
may be called the universal $m$-th order variation of Hodge structure
associated to $X$.
We will give a formula for it, together with (something equivalent to) it
Gauss-Manin
correction, in terms of $X$ itself.

First, for a sheaf or complex $A$ on a product $Y_1 \times Y_2$, we denote by
${}_i G.H^\cdot
(A), i = 1,2$, the increasing filtration on $H^\cdot(A)$ associated to the
Leray spectral
sequence
$$
E^{p,q}_2 = H^p (Y_i, R^q p_{i*} A) \Rightarrow H^\cdot (Y_1 \times Y_2, A),
$$
and set
$$
{}_i H^{r,s} (A) = \gr^r (H^{r+s} (A), {}_i G.),
$$
which may be called the generalised $r$-th Kunneth component of $H^{r+s} (A)$.
Note the
natural map
$$
\varphi : {}_1 H^{r,0} (A) \to H^r (A) \to {}_2 H^{0,r} (A) .
\tag2.3
$$
When $A = A_1 \boxtimes A_2$ this is clearly an isomorphism by Kunneth.  Note
also that
when $A$ is a complex, the stupid or Hodge filtration on $\H^\cdot (A) $
induces one
on ${}_i \H^{r,s} (A)$.
\enddemo

\proclaim{Theorem 2}  Let $X_m/R_m$ be the universal $m$-th order deformation
of a compact
K\"ahler manifold $X$ with $H^0 (\Theta_X) = 0. Then$

(i) we have Hodge filtration-preserving $R_m$-linear isomorphism
$$
\align
B_{R_m}^m H^r_{DR} (X_m/R_m) &\cong {}_1 \H^{r,0} (X <m,1>, J_m^{\cdot\cdot}
(\Theta,
\Omega_X^\cdot))\tag2.4\\
&\cong {}_2 \H^{0,r} (X <m,1>, J_m^{\cdot\cdot} (\Theta, \Omega_X^\cdot))\\
H^r_{DR} (X_m,R_m) &\simeq C_{R_m}^m ({}_1 \H^{r,0} (X<m,1> , J_m^{\cdot\cdot}
(\Theta,
\Omega^\cdot_X)) \tag2.5
\endalign
$$
(ii)  We have a Gauss-Manin isomorphisms
$$
\align
GM:  H^r_{DR} (X_m/R_m) & {\overset \sim \to\rightarrow} H^r_{DR} (X) \otimes
R_m\\
B^m H^r_{DR} (X_m/R_m) &{\overset \sim \to \rightarrow} H^r_{DR} (X) \otimes
B_0^m .
\endalign
$$
which shifts the Hodge filtration by $m$.  {\rm (We call the resulting
projection
$\overline{GM}: B^m H^r_{DR} (X_m/R_m) \to H^r_{DR} (X)$ the Gauss-Manin
projection).
\endproclaim

\demo{Proof}  As $H^r_{DR} (X_m/R_m)$ is $R_m$-free, the isomorphisms (1.1),
(1.2) show that
(2.4) and (2.5) are mutually equivalent.  By the Poincar\'e lemma we have a
$\Theta$-linear
quasi-isomorphism
$$
\C \sim \Omega_X^\cdot
$$
where $\Theta$ acts trivially on $\C$.  Hence
$$
J_m^{\cdot\cdot} (\Theta, \Omega_X^\cdot) \sim J_m^\cdot (\Theta) \boxtimes \C
\oplus \C_X
\tag2.6
$$
where $X= \Delta_X \subset X <1> \times X \subset X <m_1 1>$.  This implies
degeneration
of the Leray spectral sequences for $J_m^{\cdot\cdot} (\Theta, \Omega^\cdot)$
with respect
to both projections $p_1$ and $p_2$, and that the map $\varphi$ is (2.3)--which
clearly
respects Hodge filtrations is an isomorphism.  Also, from $p_2$ we get an
isomorphism
$$
{}_1 \H^{r,0}(J_m^{\cdot\cdot} (\Theta),\Omega^\cdot)) {\overset
\sim\to\rightarrow} H^r
(X, \C) \otimes B_0^m = H^r (X,\C) \otimes (\C \oplus \H^0 (X<m>, J_m^\cdot
(\Theta)).
$$
On the other hand in view of Theorem 1 we get via $p_2$ an isomorphism
$$
{}_2 \H^{0,r} (J_m^{\cdot\cdot} (\Theta), \Omega_X^\cdot))
{\overset\sim\to\rightarrow}
\H^r (X, B^m (\Omega^\cdot_{X_m/R_m})).
$$
As $\Omega^\cdot_{X_m/R_m}$ is $R_m$-flat, the latter may be identified with
$B^m \H^r
(X, \Omega^\cdot_{X_m/R_m}) = B^m H^r_{DR} (X_m/R_m)$, which completes the
proof of (i).
\enddemo

The proof of (ii) is based on the Cartan formula for the Lie derivative of
differential
forms
$$
L_v = i_v\circ d + d\circ i_v \qquad i_v = \ \ \text{interior
multiplication.}
$$
Let $\Omega_{X,\ \text{triv}}^\cdot$ be $\Omega_X^\cdot$ with the trivial
action of
$\Theta$ and define a map
$$
\align
M:tot (J_m^{\cdot\cdot} (\Theta, \Omega_X^\cdot)) &\to tot (J_m^{\cdot\cdot}
(\Theta, \Omega_{X, \ \text{triv}}^\cdot))\tag2.7\\
M^{i,j,k} : \lambda^j \Theta \boxtimes \Omega_X^i &\to \lambda^{j-k} \Theta
\boxtimes
\Omega^{i-k}
\endalign
$$
$$
M^{i,j,k} (v_1 \times \cdots \times v_j \times \omega) = \sum_{r_1 < \cdots <
r_k}
(-1)^{\sum{r_s}}  v_1 \times \cdots \times \hat{v}_{r_1} \cdots \times
\hat{v}_{r_k} \times \cdots \times v_j \times i_{v_{r_1}\cdots v_{r_k}}
(\omega)
$$
The Cartan formula implies that (with proper choice of signs) $Id \oplus M$ is
a morphism of
complexes, hence yields a map
$$
\align
{}_1 \H^{r,0} ( X <m,1>, J_m^{\cdot\cdot} (\Theta, \Omega_X^\cdot)) &\to {}_1
\H^{r,0}
(X<m,1>, J_m^{\cdot\cdot} (\Theta, \Omega^\cdot_{X,\ \text{triv}}))\\
&\simeq \H^r (\Omega^\cdot_X)\otimes B_0^m
\endalign
$$
which, in view of the quasi-isomorphism (2.6) is an isomorphism, and obviously
it shifts
the Hodge filtration by $m$.  Applying the $C^m$ functor then concludes the
proof.

Remark that $\Theta \oplus \Omega^\cdot_X$ has a structure of (bidifferential
bigraded) Lie
algebra and the image of $J_m^{\cdot\cdot} (\Theta, \Omega^\cdot)$ on $X<m+1>$
may
be identified with a direct summand of $J_{m+1}^{\cdot\cdot} (\Theta \oplus
\Omega^\cdot)$.
Correspondingly $B^m H^r_{DR} (X_m/R_m)$ may be identified with a direct
summand of
$\H^{0,r} (J_{m+1}^{\cdot\cdot} (\Theta \oplus \Omega^\cdot)) = T^{m+1} (H^r)$
where
$H^r$ is the cohomology bundle $H^r_{DR} (X_m/R_m)$ viewed as a geometric
vector bundle
over $\Spec (R_m)$ (cf. 1.1)
\ls
\subheading{3.  The Calabi-Yau case}
\ss
We will now apply the foregoing methods to study the period map for (the
$n$-form on)
Calabi-Yau manifolds.  So let $X$ be Calabi-Yau, i.e. $X$ is an $n$-dimensional
compact
K\"ahler manifolds with a nowhere-vanishing holomorphic $n$-form $\Phi$ (which
is then unique
up to isomorphism); we will also assume $X$ admits no holomorphic vector
fields, i.e.
$H^0 (\Theta_X) = 0$.  The pair $(X, \Phi)$ may be called a {\it measured}
Calabi-Yau
manifold (MCYM).  An isomorphism between two MCYM's $(X,\Phi)$ and $(X',\Phi')$
in an
isomorphism $f: X \to X'$ with $f^* \Phi' = \Phi$.  Thus the Lie algebra sheaf
of
infinitesimal automorphisms of a MCYM may be identified as the subsheaf
$\hat{\Theta}
\subset \Theta_X$ of {\it divergence-free} vector fields, i.e. those
annihilating
$\Phi$ via Lie derivative.  As before, we have universal $m$-th order
deformations
$X_m/R_m$ of $X$ and $(X_m/\hat{R}_m, \Phi_m)$ of $(X,\Phi)$, where $R_m$
(resp.
$\hat{R}_m$) is the algebra associated to the OS structure $\H^0 (J_m^\cdot
(\Theta))$
(resp. $\H^0(J_m^\cdot (\hat{\Theta}))$ (in fact, $\hat{R}_m$ is just the
obvious ``$m$-th order" quotient of $R_m[t]$).  The period space for $(X,\Phi)$
is just the
vector space $H = H^n_{DR} (X)$ (while for $X$ it would be the projectivisation
$\P(H))$.
The $m$-th order germ of the period map yields a ring homomorphism
$$
p_m^* : S^\cdot [H] \to \hat{R}_m.
$$
Of course $p^*_m$ is determined by its restriction $p_m^{1*}$ on $H$, hence by
the dual
$$
p_m^1 : \hat{T}^{(m)} \to H = H^*
{}.
$$
Our purpose is to give a cohomological formula for $p_m^1$.  To this end let's
replace
the DeRham complex $\Omega_X^\cdot$ by its quasi-isomorphic subcomplex
$\Omega_{X,0}^\cdot$
given by
$$
\align
\Omega^\cdot_{X,0} & = \Omega^i_X \quad i \leq n-2\\
&= \hat{\Omega}^{n-1}_X \quad i = n-1\\
&= 0 \quad i = n .
\endalign
$$
Define a map $j^\cdot_m : J^\cdot_m (\hat{\Theta}) \to \Omega^\cdot_{X,0} [n]$
by
$$
\align
j_m^k (v_1\times \cdots\times v_k) &= i_{v_1\wedge\cdots\wedge v_k} \Phi \quad
v_1 \times
\cdots \times v_k \in \lambda^k \hat{\Theta} = J_m^{-k} (\hat{\Theta}) \\
&= 0 \quad k\geq n + 1 .
\endalign
$$

\proclaim{Theorem 3}  $j_m^\cdot$ is a morphism of complexes and the associated
cohomology map
$\H^0 (j_m^\cdot): \hat{T}^m \to H$ coincides with $p_m^1$.
\endproclaim

\demo{Proof}  We begin with a geometric description of the period map $p$.
Consider the germ
$\hat{M}$ of the moduli for the MCYM $(X,\Phi)$. Over $\hat{M}$ we have the
cohomology bundle
$\Cal H = \Cal H^n$, and the Gauss-Manin isomorphism
$$
GM: \Cal H \tilde{\rightarrow} \hat{M} \times H
$$
and the Gauss-Manin projection  $\overline{GM}: \Cal H \to H$.  Now the
`tautological'
section $\Phi$ yields a cross-section $[\Phi]:\hat{M} \to \Cal H$.

In terms of these, the period map $p$ is simply given by
$$
p = \overline{GM} \cdot [\Phi] : \hat{M} \to H.
\tag3.0
$$
On $m$-th order tangent spaces we get a map
$$
p_m = T^m (p): \hat{T}^m = T^m (\hat{M}) \to T^m H = \opp\limits_1^m S^i (H).
$$
Being a morphism of OS structures, $p_m = \oplus p_m^i$ is determined by $p_m^1
= dp^m$.
As in \S1.2, we have
$$
p_m^1 = \overline{GM} \circ B^m (\Phi),
$$
where $B^m (\Phi) : \hat{T}^m \to B^m(\Cal H)$ is the map corresponding to
$\Phi$.
Clearly $B^m (\Phi)$ is the map on $\Bbb H^0$ corresponding to the morphism of
complexes
$$
u_m^\cdot : J_m^\cdot (\hat{\Theta}) \to J_m^\cdot (\Theta, \Omega^n_X)
\hookrightarrow
J_m^{\cdot\cdot} (\hat{\Theta}, \Omega^\cdot_X)
$$
given by $v_1 \times \cdots \times v_m \mapsto v_1 \times \cdots\times v_m
\times
\Phi$ (this is compatible with differentials because $L_{v_i} \Phi = 0$).  On
the other
hand $\overline{GM}$ is given by $M^{\cdot,j,j}$ as in (2.7) and clearly
$M^{n,j,j} \circ
u^j_m$, factors through $\Omega_{X,0}^j$ and coincides with $j_m^j$, hence
$p_m^1 =
\Bbb H^0 (j_m^\cdot)$ as claimed.
\enddemo

Let us now consider consequences of Theorem 3.
Denote by $\sigma^m(\Omega_{X,0}^\cdot [n])$ the $m$-th exterior
signed-symmetric tensor
power of $\Omega^\cdot_{X,0} [n]$, as complex on $X <m>$, and
set $\sigma^{m]}(\Omega_{X,0}^\cdot [n])=\opp\limits_0^m
\sigma^i(\Omega_{X,0}^\cdot [n])$
note that $j_{m}^\cdot$
naturally extends to a morphism
$$
\sigma^{m]} (j_m^\cdot): J_m^\cdot (\Theta) \to \sigma^{m]} (\Omega^\cdot_{X,0}
[n])
$$
(similarly given by interior multiplication by $\Phi$).
As $p_m$ is a morphism of OS structures, it clearly follows from Theorem 3
that:

\proclaim{Corollary 3.1}  We have $p_m = \Bbb H^0 (\sigma^{m]} (j_m^\cdot)):
\hat{T}^m
\to S^{m]} (H)$.
\endproclaim

Now let us endow $H$ with the modification of the Hodge filtration given by
setting $F_+^1 = H, F_+^\cdot = F^\cdot$ otherwise.  This filtration then
induces one
on $S^\cdot (H)$ making $p_m^*$ a filtered homomorphism inducing an isomorphism
$$
\gr^1_{F_+} S^\cdot (H) = H/F^2 H = (F^{n-1} H)^* \tilde{\rightarrow} gr^1
(\hat{R}_m) = \hat{T}^*.
$$
Note the diagram
$$
\matrix
\hat{T}^m &\to& F^m ( S^2(\hat{T}^{m-1}))\\
p_m \downarrow&&\downarrow S^2 (p_{m-1})\\
\opp\limits_1^m S^i (H)&\to& F^m (S^2 (\opp\limits_1^{m-1} S^i (H)))
\endmatrix
\tag3.1
$$
(horizontal maps being symbol or comultiplication maps).  It is essentially the
commutativity of (3.1) and its
dual that translate into equations for the image of $p_m$:  consider the
following set of
Yukawa-type elements
$$
\multline
Y_n = \left\{ a-b\cdot c \in S^\cdot (H) : a \in H, b \cdot c \in
\opp\limits_1^n S^\cdot
(H),\right. \\
\left. p_n^{1*} (a) = p_{n-1}^* (b) \cdot p_{n-1}^* (c) \right\}
\endmultline
$$

Let $S_m = S^\cdot (H)/F^{m+1} S^\cdot(H)$ be the $m$-th quotient of the
filtered ring
$S^\cdot(H)$ and $S^m = \ker (S_m \to S_{m-1})$ the $m$-th graded piece.  Let
$I_m$ (resp.
$K_m$) be the kernel of the natural surjection $S_m \to \hat{R}_m$ (resp. $S^m
\to
S^m (\hat{T}^*))$ induced by $p_m^*$.  By induction on $m\geq n$, we proceed to
define
lifts $Y_m \subseteq I_m$ of $Y_n$.  given $y\in Y_m$, start with an arbitrary
lift $y' \in S_{m+1}$.
As $p_{m+1}^* (y') \in S^{m+1} (\hat{T}^*) \subset \hat{R}_{m+1}$, there is an
uniquely
determined element
$$
z = z(y') \in S^{m+1} (gr^1_{F_+^\cdot} (H)) \subseteq S^{m+1}
\subset S_{m+1}
$$
such that $p_{m+1}^* (z) = p_{m=1}^* (y')$, and we let $Y_{m+1} = \{y' - z(y')
: y \in Y_m\}$.

Note on the other hand that $K_m$ may be easily written explicitly in terms of
the Yukawa
forms $\eta^i: S^i \hat{T} \to H^{n-i, i}$, which formally determine the map
$p^m = gr^m
p_m^* : S^m \to S^m (\hat{T}^*)$.  In fact it is clear as in the above that
$K_m$ is
generated by elements of the form $y-z$ where $y\in S^m, z \in S^m (gr^1
(H)^*),
p^m (y) = p^m (z)$.  In particular, it follows easily that for $m\geq n, K_m$
is
generated by $K_n$, i.e. $K_m = F_+^{n-m} S_m \cdot K_n$.

\proclaim{Corollary 3.2}  For all $m\geq n$ $Y_n$ constitutes a complete set of
defining
equations for the image of the $m$-th order period map, i.e. $p_m^*$ induces an
isomorphism
$$
S_m / (K_{m} + <Y_m>) \simeq \hat{R}_m.
$$
\endproclaim

\demo{Proof}  Let us denote the LHS by $A_m$.  As $p_m^*$ is compatible with
the $F_+^\cdot$
and $m$-adic filtrations, it would suffice to prove that $\gr^\cdot p_m^* :
\gr^\cdot
A_m \to \gr^\cdot \hat{R}_m = \opp\limits_0^m S^i \hat{T}$ is an isomorphism.
For
$\gr^1$ this is clear.  Consequently, it is easy to see that it would suffice
to prove
that the evident map $S^j H \to \gr^i A_m$ is surjective.  e.g. for $j =2$ a
typical
element $a \in \gr^2 A_m$ may be represented by $a_1 + a_2, a_1 \in F^2 H, a_2
\in S^2 H$.
Clearly $p_m^* (a_1) \in m^2_{\hat{R}_m}$ so can be written as $p_m^* (b) \cdot
 p_m^* (c),
b, c \in H$, so in $\gr^2 A_m$ we have $a = a_2 + bc \in \im S^2 H$.  The case
of general
$j$ follows similarly.
\enddemo

Letting $\hat{S} = \lim\limits_\leftarrow S_m  , \Hat{\Hat{R}} =
\lim\limits_\leftarrow R_m$, we conclude
$$
\hat{S} / <Y_m> \cong \Hat{\Hat{R}} .
$$
Since
the above isomorphism is homogeneous with respect to scaling the $n$-form
$\Phi$, i.e. is
compatible with $F^1 H \subset H, R_m \subset R_m$, a similar assertion can be
made about
the germ $M \subset \Bbb P (H)$.

With a slightly more analytic approach, the above construction of Schottky
relations
may be made more explicit.  Firstly, we may split the Hodge filtration, e.g.
with the
Hodge decomposition $H = \oplus H^{n-i,i}$.  Then setting
$$
H^{n-1,1}_+ = H^{n,0} \oplus
H^{n-1,1}, H^{1,n-1}_* = (H^{n-1,1}_+)^* = H^{0,n} + H^{1,n-1}
$$
we get a splitting of $F_+^\cdot$
and of the dual filtration $F_*^\cdot$ on $H = H^*$. Assigning
weight 1 (resp. $i$) to $H^{1,n-1}_*$ (resp. $H^{i,n-i}, i \geq 2$)
then induces a graded ring structure $B^\cdot$ on $B = S^\cdot (H),$ i.e.
$$
B^\cdot = \Bbb C \oplus H^{1,n-1}_* \oplus (S^2(H^{1,n-1}_*) \oplus
H^{2,n-2}) \oplus ...
$$
Note that the graded subring $B^{m]}_1 = S^{\leq m}(H^{1,n-1}_*) \subseteq
B^{m]}
$ is mapped by $p_m^*$ isomorphically onto $\hat{R}_m$, and that $B_1^\cdot$
may
be viewed as part of the $\bar{\partial}$ -cohomology of $S^\cdot _{\Bbb
C}(A_X^
{\cdot ,\cdot}).$
$$
$$
Now identifying $B^{m]} = \hat{R}_m,$ the associated graded $gr^i(p^*_m)$
(which is
independent of $m\geq i$) is given by the 'dual Yukawa map' $\eta ^{i*}$, i.e.
corresponds to a map given on the form level by
$$
\eta ^{i*} : A^{i,n-i} \simeq A^{n-i,i *} \lra S^i (A^{n-1,1})^* \simeq
S^i (A^{1,n-1}), i\geq 2
$$
$$
\eta^{1*} = \text{id}
$$
(recall that $\eta ^i$ itself is given on the form level by
$$
\eta^i (\gamma_1...\gamma_i) = \text{skew-symmetrisation of }i_{\Phi ^{-i+1}}(
\gamma_1 \otimes ... \otimes \gamma_i)),
$$
$\Phi ^{-1} \in H^0 (\Lambda ^n (\Theta ))$ being the dual of $\Phi$; $\eta
^{i*}$
on the other hand appears to be a somewhat 'deeper' object and doesn't seem to
be given by such a simple purely local expression).
$$
$$
Now to construct Schottky relations, i.e. a basis for $I_m = \text{ker}(S^{m]}
\lra \hat{R}_m )$, we proceed inductively. As $K_m = I_m \cap S^m$
may be assumed known, it would suffice to describe a canonical lift of any
$y\in I_m$ to $I_{m+1}$. We may assume $y$ has no component in $H^{0,n}$. On
the
form level we may write
$$
\eta^* (y) = \bar{\partial}(f) = \sum\limits _{\alpha} \bar{\partial}f_1^\alpha
\otimes f_2^\alpha
 \otimes ...\otimes f_{m_\alpha}^\alpha,
$$
where $f^\alpha _1 \in A^{1,n-1}, f^\alpha _j \in A^{1,n-1}, j\geq 2, 1\leq
m_\alpha
\leq m,$ and the expression $f = \Sigma \Pi f^\alpha _i$ may be written in term
of $y$
via a suitable Green's operator. We may then define a polynomial $q_y = \sum
q_y^\alpha$
on $\hat T$ by the formula
$$
q^\alpha _y (u_0,...,u_{m_\alpha +1}) = \text{u-symmetrisation of} \ \ \
(\int_X [u_0,u_1]f^\alpha
 _1)(\int_X u_2f^\alpha _2)...(\int_X u_{m_\alpha }f^\alpha _{m_\alpha})
$$
where $u_i \in \hat {T}$,
$$
[u_i,u_j] = \partial\eta (u_i,u_j)\in A^{n-1,2}
$$
($[.,. ] $ being the
map induced by the Lie bracket
, and the latter equality being a form of the 'Tian-Todorov Lemma'). Then $q_y$
corresponds to an element $z_y\in B_1^{m+1]}$
and we set
the lift of $y$ as
$$
y' = y-z_y \in I_{m+1}.
$$
\vskip1in
{\it Remarks}  1. For $n=2$, the construction may be simplified and refined.
Indeed it easy
to write down a subcomplex $K^{\cdot} \subseteq \sigma ^{\cdot} (\Omega
_{X,0}^{\cdot}[n])$
which contains- and indeed coincides with- the image of $\sigma
^{\cdot}(j^{\cdot})$:
e.g. $K^{-2} = \left \{ (a,b \times c) \in \Cal O \times \lambda ^2 (\Omega
^1_{X,0}):
a \Phi = b \wedge c \right \}$, etc. and we have
$$
\Bbb H ^0 (K^{\geq -m})^* = S^{\cdot}(H)/(\Phi - q),
$$
where $q \in S^2(H)$ is the canonical quadratic form . Clearly then $p^*_m$
factors through
an isomorphism
$$
\Bbb H ^0 (K^{\geq -m})^* \simeq \hat R_m.
$$
Thus in this case there is a single quadratic Schottky relation $\Phi -q$, i.e.
the familiar 'period quadric' of K3 surface theory.
\ss
2. Identifying $\hat R_m = B^{m]}_1$, the period map $B^{m]} \lra B^{m]}_1$
fails to
respect the respective gradings, and it is precisely this failure that, for
$n=3$,
is supposed to be related to the quantum cohomology of the mirror manifold.

\vskip2in

\subheading{References}
\ss

\item{[D]}  Deligne, P.:  Th\'eor\'eme de Lefschetz et crit\`eres de
d\'eg\'en\'erescence de
suites spectrales'.  Publ. Math. IHES {\bf 35} (1968), 197--226.
\ss
\item{[G]}  Griffiths, P.:  `Topics in transcendental algebraic geometry' Ann.
of Math. Studies.
\ss
\item{[R1]}  Ran, Z.:  `Canonical infinitesimal deformations' (preprint).
\ss
\item{[R2]} $\underline{\hskip.5in}$  `On the local geometry of moduli spaces
of vector
bundles (preprint).

\bye